\documentclass{amsart}

\usepackage{amsmath, amssymb}
\usepackage{amscd}
\usepackage{verbatim}
\usepackage{tikz}
\usetikzlibrary{arrows,chains,matrix,positioning,scopes}

\makeatletter
\tikzset{join/.code=\tikzset{after node path={%
\ifx\tikzchainprevious\pgfutil@empty\else(\tikzchainprevious)%
edge[every join]#1(\tikzchaincurrent)\fi}}}
\makeatother

\tikzset{>=stealth',every on chain/.append style={join},
         every join/.style={->}}

\newtheorem{Theorem}{\sc Theorem}[section]
\newtheorem{Lemma}[Theorem]{\sc Lemma}
\newtheorem{Proposition}[Theorem]{\sc Proposition}
\newtheorem{Corollary}[Theorem]{\sc Corollary}
\newtheorem{Definition}[Theorem]{\sc Definition}
\newtheorem{Example}[Theorem]{\sc Example}
\newtheorem{Remark}[Theorem]{\sc Remark}
\newtheorem{Note}[Theorem]{\sc Note}
\newtheorem{Question}[Theorem]{\sc Question}
\newtheorem{ass}[Theorem]{\sc Assumption}
\newcommand{\bt}{\begin{Theorem}}
\def\beginlem{\begin{Lemma}}
\def\beginprop{\begin{Proposition}}
\def\begincor{\begin{Corollary}}
\def\begindef{\begin{Definition}}
\def\beginexamp{\begin{Example}}
\def\beginrem{\begin{Remark}}
\def\beginq{\begin{Question}}
\def\beginass{\begin{ass}}
\def\beginnote{\begin{Note}}
\newcommand{\et}{\end{Theorem}}
\def\endlem{\end{Lemma}}
\def\endprop{\end{Proposition}}
\def\endcor{\end{Corollary}}
\def\enddef{\end{Definition}}
\def\endexamp{\end{Example}}
\def\endrem{\end{Remark}}
\def\endq{\end{Question}}
\def\endass{\end{ass}}
\def\endnote{\end{Note}}

\def\textmatrix#1&#2\\#3&#4\\{\bigl({#1 \atop #3}\ {#2 \atop #4}\bigr)}
\def\dispmatrix#1&#2\\#3&#4\\{\left({#1 \atop #3}\ {#2 \atop #4}\right)}

\newcommand{\la}{\langle}
\newcommand{\ra}{\rangle}

\newcommand{\Bd}{\mathbb{B}_d}
\newcommand{\clU}{\mathcal{U}}
\newcommand{\clY}{\mathcal{Y}}
\newcommand{\clV}{\mathcal{V}}
\newcommand{\clW}{\mathcal{W}}
\newcommand{\ot}{\otimes}
\newcommand{\vp}{\varphi}
\newcommand{\op}{\oplus}
\newcommand{\sduy}{\mathcal{S}_d(\mathcal{U},\mathcal{Y})}
\newcommand{\sdvw}{\mathcal{S}_d(\mathcal{V},\mathcal{W})}
\newcommand{\al}{\alpha}
\newcommand{\be}{\beta}
\newcommand{\hkp}{H(k_\vp)}
\newcommand{\hks}{H(k_\psi)}
\newcommand{\hdty}{H^2_d(\mathcal{Y})}
\newcommand{\hdtw}{H^2_d(\mathcal{W})}
\newcommand{\itb}{I_{H^2_d} \otimes \beta}
\def ~{\hspace{1mm}}
\newcommand{\Dp}{\mathcal{D}_\varphi}
\newcommand{\Ds}{\mathcal{D}_\psi}
\newcommand{\Pd}{\Gamma^d}
\newcommand{\clM}{\mathcal{M}}
\newcommand{\vpj}{\varphi_J}
\newcommand{\sjp}{\psi_{J^\prime}}
\newcommand{\clH}{\mathcal{H}}
\newcommand{\clN}{\mathcal{N}}

\begin{document}

\title[COINCIDENCE OF SCHUR MULTIPLIERS]{COINCIDENCE OF SCHUR MULTIPLIERS OF THE DRURY-ARVESON SPACE}

\author[A.~Bhattacharya]{Angshuman Bhattacharya}
\author[T.~Bhattacharyya]{Tirthankar Bhattacharyya}

\address{Department of Mathematics, Indian Institute of Science, Bangalore-560012, India}
\email{angshu@math.iisc.ernet.in, tirtha@member.ams.org}


\begin{abstract}
In a purely multi-variable setting (i.e., the issues discussed in
this note are not interesting in the single variable operator
theory setting), we show that the coincidence of two operator
valued Schur class multipliers of a certain kind on the
Drury-Arveson space is characterized by the fact that the
associated colligations (or a variant, obtained canonically) are
`unitarily coincident' in a sense to be made precise in this
article.
\end{abstract}

\maketitle

\section{introduction}

Given a Hilbert space $\clU$ (all Hilbert spaces in this note are over the complex field and are separable),
$H^2_d(\clU)$  denotes the reproducing kernel Hilbert space corresponding to the kernel
$$\hspace{15mm}k(z,w)=\frac{I_\clU}{1-\la z,w \ra} \hspace{5mm}; \hspace{5mm}(z,w)\in\Bd \times \Bd.$$
It consists of $\clU$ valued holomorphic functions on the Euclidean unit ball $\mathbb{B}_d$ in $\mathbb{C}^d$.

\vspace{5mm}

   If $\clY$ is another Hilbert space, then a multiplier $\vp$ is a $\mathcal{B}(\clU,\clY)$
   valued holomorphic function on $\Bd$ such that $\vp f \in H^2_d (\clY)$ \hspace{1mm} for all $f\in H^2_d (\clU)$.
   Here $\vp f$ denotes the pointwise product. It is a well known consequence of the closed graph theorem that the linear map
   $f \mapsto \vp f$ is a bounded operator $M_\vp$ from $H^2_d(\clU)$ into $H^2_d(\clY) $.
   The multipliers form a Banach space with the norm
   $$\| \vp\|:=\|M_\vp \|.$$
   The closed unit ball of this Banach space is called the Schur class $\sduy$. The following theorem to be found in \cite{AMCcp} and \cite{AMc} describes the Schur class.

\bt

Let $\vp$ be a $\mathcal{B}(\clU,\clY)$ valued function defined on $\Bd$. Then the following are equivalent :
\begin{enumerate}
  \item $\vp \in \sduy$.
  \item The kernel $k_\vp  : \Bd \times \Bd \rightarrow \mathcal{B}(\clY)$ given by $$k_\vp (z,w)= \frac{I_\clY - \vp(z) \vp(w)^*}{1-\la z,w \ra}$$ is positive semi-definite.
  \item There exists an auxiliary Hilbert space $\mathcal{X}$ and a unitary connecting operator (or colligation) $U$ of the form
  $$U=\left(
        \begin{array}{cc}
          A & B \\
          C & D \\
        \end{array}
      \right) = \left(
                  \begin{array}{cc}
                    A_1 & B_1 \\
                    . & . \\
                    . & . \\
                    A_d & B_d \\
                    C & D \\
                  \end{array}
                \right) : \left(
                            \begin{array}{c}
                              \mathcal{X} \\
                              \clU \\
                            \end{array}
                          \right) \longrightarrow \left(
                                                    \begin{array}{c}
                                                      \mathcal{X}^d \\
                                                      \clY \\
                                                    \end{array}
                                                  \right)
  $$ such that $\vp(z)$ can be realized as $\vp(z)=D+C(I_{\mathcal{X}} - ZA)^{-1} ZB$, where $Z$ denotes the row tuple $(z_1 I_{\mathcal{X}},...,z_d I_{\mathcal{X}})$ for $z\in\Bd$.
  \item There exists an auxiliary Hilbert space $\mathcal{X}$ and a contractive colligation operator $U$ such that $\vp$ can be realized in the same form as above.
\end{enumerate}

\et

A realization, as in (3) above, of a Schur multiplier is commutative if the operators $A_i$ commute among themselves. Not all Schur multipliers admit commutative realizations. See Proposition 3.3 of \cite{BBF2}.

\begindef

The colligation operators
$$ U_1 = \left(
          \begin{array}{cc}
            A_1 & B_1 \\
            C_1 & D_1 \\
          \end{array}
        \right) : \left(
                    \begin{array}{c}
                      K_1 \\
                      \clU \\
                    \end{array}
                  \right) \rightarrow \left(
                                        \begin{array}{c}
                                          K_1^d \\
                                          \clY \\
                                        \end{array}
                                      \right) \\
                                      \text{ and } U_2 = \left(\begin{array}{cc}
            A_2 & B_2 \\
            C_2 & D_2 \\
          \end{array}
        \right) : \left(
                    \begin{array}{c}
                      K_2 \\
                      \clV \\
                    \end{array}
                  \right) \rightarrow \left(
                                        \begin{array}{c}
                                          K_2^d \\
                                          \clW \\
                                        \end{array}
                                      \right)$$

are said to be `unitarily coincident' if there exist unitary operators
$$\Lambda : K_1 \rightarrow K_2, \;\; \Omega_1 : \clU \rightarrow \clV, \mbox{ and } \Omega_2 : \clY \rightarrow \clW$$
satisfying
$$\left(
    \begin{array}{cc}
      \Lambda^d & 0 \\
      0 & \Omega_2 \\
    \end{array}
  \right) \left(
            \begin{array}{cc}
              A_1 & B_1 \\
              C_1 & D_1 \\
            \end{array}
          \right) = \left(
                      \begin{array}{cc}
                        A_2 & B_2 \\
                        C_2 & D_2 \\
                      \end{array}
                    \right) \left(
                              \begin{array}{cc}
                                \Lambda & 0 \\
                                0 & \Omega_1 \\
                              \end{array}
                            \right).
$$
\enddef

 The reproducing kernel Hilbert space $H(k_\vp)$ corresponding to a $\vp \in \sduy$ is known as the de Branges-Rovnyak space.

\vspace{2mm}

Unitary coincidence of two colligations is a concept which brings
together many facets of the realization formulae for $\vp$ and
$\psi$. It shows that the de Branges-Rovynak spaces for $\vp$ and
$\psi$ are isomorphic as Hilbert modules. In case $\vp$ and $\psi$
happen to be characteristic functions of commuting contractive
tuples, then this is all that needs to be said, see \cite{BES2}.
But when they are not, then there is an interesting extra
ingredient $B^\vp$ (respectively $B^\psi$) as defined in \cite{BBF1}. The beauty of the
following analysis lies in the fact that the coincidence of $\vp$
and $\psi$ implies unitary equivalence of $B^\vp$ and $B^\psi$
with a suitable modification if needed. The kernel spaces
$\clU^0_\vp$ and $\clV^0_\psi$ (defined in Section 2) play crucial roles in influencing
this unitary equivalence. To what extent this analysis is
dependent on the domain is an interesting question. Realization
formulae have been proved in much more generality, see \cite{DM}.

\vspace{2mm}

The space $H(k_\vp)$ is contractively included in $H^2_d(\clY)$. In general, the de Branges-Rovnyak space $H(k_\vp)$, as an algebraic subspace of $H^2_d(\clY)$ may not be invariant under the canonical backward shift tuple $M_z^* = (M_{z_1}^* ,...,M_{z_d}^* )$ on $H^2_d(\clY)$.

\textit{In this note we consider only those Schur multipliers} $\vp$ \textit{for which} $H(k_\vp)$ is $M_z^*$ \textit{invariant and the following difference quotient inequality holds for all} $f\in H(k_\vp)$ :
$$\sum_{j=1}^d \|(M_{z_j}^* \ot I_\clY)f\|_{H(k_\vp)}^2 \leq \|f\|^2_{H(k_\vp)} - \|f(0)\|^2_\clY.$$
See \cite{BBF2} for more on such $\vp$.

\vspace{2mm}

In single variable theory, the de Branges-Rovnyak space is always invariant under the backward shift operator and satisfies the one variable analog of the difference quotient inequality. To consider Schur functions satisfying the multivariable inequality is a natural extension of the one variable theory.

\begin{Definition} \label{functional-model-realization}

If a realization $\textmatrix A&B\\C&D\\$
of a Schur class function $\varphi$ is such that
$A_j = M_{z_j}^*|_{H(k_\vp)}$ for $j=1,...,d$, $Cf = f(0)~ \text{for}~ f\in H(k_\vp) $ and $D = \vp(0)$
then it is called a functional model realization. \end{Definition}

Let us note that the Schur functions admitting functional model realizations are precisely the ones which are considered in this paper as per the assumptions before, \cite{BBF2}.

The difference quotient inequality is then equivalently expressed as $$\sum_{j=1}^d A_j^* A_j + C^* C \leq I_{H(k_\vp)}.$$
A functional model realization of a Schur multiplier is not unique.

\section{Preliminaries and notations}

\subsection{Coincidence} We define the notion of coincidence of two Schur multipliers.

\begindef
 A $\vp \in \sduy$ is said to \textbf{coincide} with $\psi \in \sdvw$, if
there exist unitary operators $\al :\clU \rightarrow \clV$ and
$\be : \clY \rightarrow \clW$ such that $\be \vp (z) = \psi (z)
\al$ holds for all $z\in\Bd$.
\enddef

In this section, $\vp \in \sduy$ and $\psi \in \sdvw$ are
two such functions which coincide.

\beginlem

There is a unique linear isomorphism $\Gamma : \hkp
\to \hks$ such that
$$\Gamma(\sum_{t=1}^m k_\vp (\cdot,w_t)y_t ) =
\sum_{t=1}^m k_\psi (\cdot,w_t) \be y_t .$$

Moreover, if we identify the Hilbert spaces $\hdty$ and $H^2_d
\otimes \clY$ (and similarly for $\hdtw$), then
$\Gamma=(\itb)|_{\hkp}$.

\endlem

Proof : Uniqueness is clear because of density of the
vectors
$$ \sum_{t=1}^m k_\vp (\cdot,w_t)y_t.$$
For existence we need to show that
$$
\left\| \sum_{t=1}^m k_\psi (\cdot,w_t) \be y_t \right \| ^2_{\hks} =\left \| \sum_{t=1}^m k_\vp (\cdot,w_t)y_t \right\|^2_{\hkp}$$
which easily follows from coincidence. The rest is straightforward
computation.

To prove that $\Gamma=(\itb)|_{\hkp}$, note that, as a vector
$$k_\vp (\cdot,w)=\frac{I_\clY - \vp(\cdot) \vp(w) ^*}{1-\la
\cdot, w \ra} y \hspace{2mm} \in \hspace{1mm} \hdty.$$

Thus, $$(\itb)\hspace{2mm}k_\vp (\cdot,w)y
\hspace{2mm}\in\hspace{1mm} \hdtw.$$ Let $h\in \clW$ and $k
(\cdot,w^\prime)h$ denote an elementary tensor in $\hdtw$. Then we
have,

\begin{eqnarray*}
\la~ (\itb)k_\vp (\cdot,w)y~,~k (\cdot,w^\prime)h~\ra_{\hdtw}&=&\la~ k_\vp (\cdot,w) y~,~ (\itb)^* k (\cdot,w^\prime)h~ \ra_{\hdty} \\
&=& \la~ k_\vp (\cdot,w)y~,~k (\cdot,w^\prime)\be^*y~ \ra_{\hdty} \\
&=& \la~ k_\vp(w^\prime,w)y~,~\be^*h ~\ra_\clY \\
&=& \la~ \be k_\vp(w^\prime,w)y~,~h \ra_\clW \\
&=& \la~ k_\psi(w^\prime,w)\be y~,~h \ra_\clW\\
&=& \la~ k_\psi(\cdot,w)\be y~,~k (\cdot,w^\prime)h \ra_{\hdtw}~.
\end{eqnarray*}

~

The second equality from the end in the above follows from the
coincidence of the two Schur multipliers. This shows that
$\Gamma=(\itb)|_{\hkp}$ on a dense subspace of $\hkp$. Finally, by
contractive inclusion of $\hks$ in $\hdtw$ one has the required
result. \qed

\subsection{The two de Branges - Rovnyak spaces} Now the following intertwining relation is immediate:
$$\Gamma ( M_{z_j}^* \ot I_\clY )|_{H(k_\vp)} = ( M_{z_j}^* \ot I_\clW
)|_{H(k_\psi)} \Gamma .$$
From the above, it is easy to see that if
one of the Schur multipliers is such that its corresponding de
Branges-Rovnyak space is backward shift invariant, the elements of
which satisfy the difference quotient inequality then the same is
also true for the other Schur multiplier which is coincident to
it.

Consider the Hilbert space $\hkp^d$. Let $W^* k_\vp (\cdot,
w)y :=\left(
                               \begin{array}{c}
                                 \overline{w}_1 k_\vp (\cdot, w)y \\
                                 . \\
                                 . \\
                                 \overline{w}_d k_\vp (\cdot, w)y \\
                               \end{array}
                             \right) \in \hkp^d$ and let $ \Dp=\overline{\mathrm{span}}~ \{~W^* k_\vp (\cdot, w)y~;~w\in\Bd, y\in\clY~ \}.$
  Similarly we define $\Ds$.
On $\hkp^d$, let $$ \tilde{\Gamma} = \left(
                           \begin{array}{ccc}
                             \Gamma & \cdots & 0 \\
                              & \ddots &  \\
                             0 & \cdots & \Gamma \\
                           \end{array}
                         \right) \mbox{ and } A^\vp = \left(
                                                     \begin{array}{c}
                                                       A^\vp_1 \\
                                                       \vdots \\
                                                       A^\vp_d \\
                                                     \end{array}
                                                   \right).$$
By the action of the unitary $\Gamma$, one has $\tilde{\Gamma}
(\Dp)=\Ds$ and thus $\tilde{\Gamma} (\Dp^\perp)=\Ds^\perp$. Let
$R^\vp : \Dp \rightarrow \clU$ be given by $W^* k_\vp (\cdot, w)y
\longmapsto (\vp(w)^* - \vp(0)^*)y .$ Next, we define
some operators considered in Section 2 of \cite{BBF1} and follow
the notations therein closely for our convenience. In what follows,
$\textmatrix A&B\\C&D \\$ is a functional model realization as in
Definition \ref{functional-model-realization}.

\begin{itemize}
  \item $T_{11}^\vp : \Dp^\perp \rightarrow \hkp$ is given by
      $T_{11}^\vp = A^*|_{\Dp^\perp}$, where $A^*=({A^*}_1,...,{A^*}_d)$.
  \item $T_{12}^\vp : \Dp \op \clY \rightarrow \hkp$ is given
      by $T_{12}^\vp = \left(
                                                                         \begin{array}{cc}
                                                                           {A^*}|_{\Dp} & C^ * \\
                                                                         \end{array}
                                                                       \right)
  $.
  \item $T_{22}^\vp : \Dp \op \clY \rightarrow \clU$ is given
      by $T_{22}^\vp = \left(
                                                                         \begin{array}{cc}
                                                                           R^\vp & \vp(0)^* \\
                                                                         \end{array}
                                                                       \right)
  $.
\end{itemize}

Note that $\overline{\mathrm{ran}}~ T^\vp_{22} \subset  \clU^{0
\perp}_\vp$. We know that any functional model realization of
$\vp$ is the adjoint of the colligation written in
terms of the above operators :
$$\left(
    \begin{array}{cc}
      T_{11}^\vp & T_{12}^\vp \\
      X & T_{22}^\vp \\
    \end{array}
  \right) : \left(
              \begin{array}{c}
                \Dp^\perp \\
                \Dp \\
                \clY \\
              \end{array}
            \right) \rightarrow \left(
                                  \begin{array}{c}
                                    \hkp \\
                                    \clU \\
                                  \end{array}
                                \right)
$$
where $\hkp^d = \Dp^\perp \op \Dp$ and $X : \Dp^\perp \rightarrow
\clU$ is the non-unique component of the operator $B^\vp = \left(
                                                                                                                               \begin{array}{cc}
                                                                                                                                 X & R^\vp \\
                                                                                                                               \end{array}
                                                                                                                             \right) ^*
 : \clU \rightarrow \hkp^d$ appearing in the colligation. See Section 2 of \cite{BBF1} for more details.
Similarly, one may consider the analogous operators for $\psi \in
\sdvw$. We record a few simple intertwining properties of the
operators $T^\vp_{ij}$ and $T^\psi_{ij}$ for $i,j=1,2$. They are
easy consequences of the intertwining property of $\Gamma$ and coincidence of the Schur
multipliers.

$$ \Gamma~T_{11}^\vp =  T_{11}^\vp ~
      \tilde{\Gamma}|_{\Dp^\perp} \vspace{0.2cm}, \;\; \Gamma~T_{12}^\vp = T_{12}^\psi~\left(
                                             \begin{array}{cc}
                                             \tilde{\Gamma}|_{\Dp} & 0 \\
                                             0 & \be
                                           \end{array}
                                           \right),\;\;
\al~T_{22}^\vp = T_{22}^\psi~\left(
                                          \begin{array}{cc}
                                            \tilde{\Gamma}|_{\Dp} & 0 \\
                                            0 & \be \\
                                          \end{array}
                                        \right).$$

Define \begin{eqnarray}
  G^\vp_1 : \overline{\mathrm{ran}}~(I_{\hkp}
      - T_{12}^\vp T_{12}^{\vp
      *})^{\frac{1}{2}} &\rightarrow&  \Dp^\perp \\
      (I_{\hkp} - T_{12}^\vp T_{12}^{\vp *})^{\frac{1}{2}} f
      &\longmapsto&  T_{11}^{\vp *} f ,~~~~ f \in
      \hkp.
      \end{eqnarray}
  and

\begin{eqnarray}
G^\vp_2 : \overline{\mathrm{ran}}~(I_{\Dp \op \clY} -
      T_{12}^{\vp *}  T_{12}^\vp)^{\frac{1}{2}} &\rightarrow&
      \clU \\
      (I_{\Dp \op \clY} - T_{12}^{\vp *}
      T_{12}^\vp)^{\frac{1}{2}} l &\longmapsto& T_{22}^\vp l,
      ~~~~l \in \Dp \op \clY.
      \end{eqnarray}

Note that $\overline{\mathrm{ran}}~G^\vp_2 = \overline{\mathrm{ran}}~ T^\vp_{22}
\subset \clU^{0 \perp}_\vp$. Similarly one has the operators
$G^\psi_1$ and $G^\psi_2$.

\beginlem

There are intertwining relations of $G^\vp_j$ and $G^\psi_j$ for
$j=1,2$ as given below :

$$ (~\Pd|_{\Dp^\perp})~G^\vp_1 = G^\psi_1~\Gamma \mbox{ and } \al ~ G^\vp_2 = G^\psi_2~\left(
                                    \begin{array}{cc}
                                      \Pd|_{\Dp} & 0 \\
                                      0 & \be \\
                                    \end{array}
                                  \right) .$$
\endlem

Proof: The proof follows from the following relations

\begin{itemize}
  \item $\Gamma~ (I_{\hkp} - T_{12}^\vp T_{12}^{\vp *}) =
      (I_{\hks} - T_{12}^\psi T_{12}^{\psi *})~ \Gamma$
      \vspace{0.5cm}
  \item $\left(
           \begin{array}{cc}
             \Pd|_{\Dp} & 0 \\
             0 & \be \\
           \end{array}
         \right)~(I_{\Dp \op \clY} - T_{12}^{\vp *}
         T_{12}^\vp) = (I_{\Ds \op \clW} - T_{12}^{\psi *}
         T_{12}^\psi)~ \left(
           \begin{array}{cc}
             \Pd|_{\Dp} & 0 \\
             0 & \be \\
           \end{array}
         \right) $
\end{itemize} \qed

\subsection{Some new colligations.} Let $\clU^0_\vp = \{~u\in\clU~ |~\vp(z)u \equiv 0~\}$. $\clV^0_\psi$ is similarly defined. Note that
by virtue of coincidence one has $\al (~\clU^0_\vp) = \clV^0_\psi$
and therefore $\al (~\clU^{0 \perp }_\vp) = \clV^{0 \perp}_\psi$.

Let $\vp \in \sduy$ be such that it admits coisometric functional
model realizations. Let $U^\vp = \left(
                                                                                                     \begin{array}{cc}
                                                                                                       A^\vp & B^\vp \\
                                                                                                       C^\vp & \vp(0) \\
                                                                                                     \end{array}
                                                                                                   \right)
$ be a coisometric functional model colligation of $\vp$. Then the
operator
$$U^{\vp *} = \left(
          \begin{array}{cc}
            {A_\vp^*} & C^{\vp *} \\
            B^{\vp *} & \vp(0)^* \\
          \end{array}
         \right) : \left(
                     \begin{array}{c}
                       \hkp^d \\
                       \clY \\
                     \end{array}
                   \right) \rightarrow \left(
                                         \begin{array}{c}
                                           \hkp \\
                                           \clU \\
                                         \end{array}
                                       \right)
         $$is an isometric operator. Now, for any auxiliary Hilbert space $\clM$, the colligation operator
$$\left(
    \begin{array}{ccc}
      A^\vp & B^\vp & 0 \\
      C^\vp & \vp(0) & 0 \\
    \end{array}
  \right) : \left(
              \begin{array}{c}
                \hkp \\
                \clU \\
                \clM \\
              \end{array}
            \right)
   \rightarrow \left(
                                  \begin{array}{c}
                                    \hkp^d \\
                                    \clY \\
                                  \end{array}
                                \right)$$
                                is also coisometric. The transfer function of the colligation above is given by $\vp_\clM (z) =\left[
                                                              \begin{array}{cc}
                                                                \vp(z) & 0 \\
                                                              \end{array}
                                                            \right]
$. Obviously, one has $\vp_\clM (z) \in \mathcal{S}_d (\clU \op
\clM , \clY)$.

The colligation above is a coisometric functional model
realization of $\vp_\clM (z)$. To see this, one simply notes that
$$k_{\vp_{\clM}} (z,w) = \frac{I_\clY - \left[
                                          \begin{array}{cc}
                                            \vp(z) & 0 \\
                                          \end{array}
                                        \right] \left[
                                                  \begin{array}{c}
                                                    \vp(w)^* \\
                                                    0 \\
                                                  \end{array}
                                                \right]
 }{1-\la z,w\ra} = \frac{I_\clY - \vp(z) \vp(w)^*}{1-\la z,w\ra} = k_\vp (z,w)~.$$
Thus, $H(k_{\vp_\clM}) = \hkp$ and hence the claim. We denote this
colligation by $U^\vp_\clM$.

It is easy to see that $(\clU \op \clM)^0_{\vp_{\clM}} =
\clU^0_\vp \op \clM$.

 Let $\clN$
be the Hilbert space
$$\clN = \clU \op \overline{\mathrm{ran}}~(I_{\Dp^\perp} - G_1^{\vp} G_1^{\vp *} )^{\frac{1}{2}}~.$$
Consider the projection  map $J : \clN \rightarrow \clU^{0
\perp}_\vp$ as a partial isometry. Its adjoint $J^* : \clU^{0
\perp}_\vp \rightarrow \clN$ is the inclusion map. The
$\mathcal{B}(\clN,\clY)$-valued function $\vp_J (z) := \vp(z) J$, on $\Bd$ and the operator valued kernel
$$k_{\vp_J} (z,w) = \frac{I_\clY - \vp_J (z) \vp_J (w)^*}{1- \la z,w \ra} $$
on $\Bd \times \Bd$ satisfy
$$k_{\vp_J} (z,w) = \frac{I_\clY - \vp_J (z) \vp_J (w)^*}{1- \la z,w \ra} = \frac{I_\clY - \vp(z)JJ^*\vp(w)^*}{1-\la z,w \ra} =
\frac{I_\clY - \vp(z)\vp(w)^*}{1-\la z,w \ra} = k_\vp (z,w)~.$$

Thus $\vp_J \in \mathcal{S}_d(\clN,\clY)$. As a consequence we
have $H(k_{\vp_J}) = \hkp$ and therefore $\vpj$ also admits
functional model realizations. Also note that
$$\clN^0_{\vp_J} =
\clU^0_\vp \op \overline{\mathrm{ran}}~(I_{\Dp^\perp} - G_1^{\vp} G_1^{\vp
*} )^{\frac{1}{2}}.$$

The following are easy consequences of the definition of $J$:

$$T_{11}^{\vpj} = T_{11}^\vp;  \;
T_{12}^{\vpj} = T_{12}^\vp; \;
T_{22}^{\vpj} = J^* T_{22}^{\vp}; \;
G_1^{\vpj} = G_1^{\vp}; \;
G_2^{\vpj}=J^* G_2^{\vp}.$$

Let $U^\vp$ be a fixed functional model realization of $\vp$. We
know that a functional model realization is always weakly
coisometric. Then from Theorem 2.7 of \cite{BBF1} we know that the
realization has an explicit description and the adjoint of the
non-unique operator $B^\vp$ is given by
$$\left(
    \begin{array}{cc}
      -G_2^\vp T_{12}^{\vp *} G_1^{\vp *} + \xi^\vp (I_{\Dp^\perp} - G_1^{\vp} G_1^{\vp *} )^{\frac{1}{2}} & R^\vp \\
    \end{array}
  \right) : \Dp^\perp \op \Dp \rightarrow \clU,
$$
for some contraction $\xi^\vp : \overline{\mathrm{ran}}~(I_{\Dp^\perp} -
G_1^{\vp} G_1^{\vp *} )^{\frac{1}{2}} \rightarrow \clU^0_\vp$. Let
$\xi^\vp_{iso}$ denote the isometric operator defined by
$$\xi^\vp_{iso} : \overline{\mathrm{ran}}~(I_{\Dp^\perp} - G_1^{\vp} G_1^{\vp *} )^{\frac{1}{2}} \rightarrow \clU^0_\vp \op
\overline{\mathrm{ran}}~(I_{\Dp^\perp} - G_1^{\vp} G_1^{\vp *} )^{\frac{1}{2}}$$
$$h \longmapsto \xi^\vp h \op \Delta_{\xi^\vp} h~.$$

 By using the sufficiency
criterion of Theorem 2.7 (2) of \cite{BBF1}, we construct the
following coisometric functional model realization of $\vpj$ :

$$U^{\vpj *} = \left(
               \begin{array}{cc}
                 {A_\vp^*} & C^{\vp *} \\
                 B^{\vpj *} & J^* \vp(0)^* \\
               \end{array}
             \right) : \left(
                         \begin{array}{c}
                           \hkp^d \\
                           \clY \\
                         \end{array}
                       \right) \rightarrow \left(
                                             \begin{array}{c}
                                               \hkp \\
                                               \clN \\
                                             \end{array}
                                           \right)~,
$$
where,
$$B^{\vpj *} = \left(
                      \begin{array}{cc}
                        [-J^*G_2^\vp T_{12}^* G_1^{\vp *} + \xi_{iso}^\vp (I_{\Dp^\perp} - G_1^{\vp} G_1^{\vp *} )^{\frac{1}{2}} ] \hspace{4mm} & J^* R^\vp \\
                      \end{array}
                    \right) : \Dp^\perp \op \Dp \rightarrow \clN~.
$$

\section{The main result and its proof}

The main result of this note is the following theorem.

\bt

Let $\vp\in\sduy$ and $\psi\in\sdvw$ be such that both admit functional model realizations. Then the following are equivalent :
\begin{enumerate}
\item $\vp$ coincides with $\psi$.
\item $\mathrm{dim}~\clU^0_\vp = \mathrm{dim}~\clV^0_\psi$ and for any two arbitrary choices of functional model colligations
$U^\vp$ and $U^\psi$ of $\vp$ and $\psi$ respectively, there exists an auxiliary Hilbert space $\clH$
such that either the canonically obtained realization $U^{\vpj}_{\clH}$ is unitarily coincident to
$U^{\sjp}$ or $U^{\sjp}_{\clH}$ is unitarily coincident to $U^{\vpj}$. Here $U^{(\cdot)}_{\clH}$ is defined as before.
\end{enumerate}

\et

As usual, one of the implications is easier than the other, and
that is $(2) \Rightarrow (1)$. Let condition $(2)$ hold. Without loss of generality let us assume that $U^{\vpj}_{\clH}$ is
unitarily coincident to $U^{\sjp}$. Then their corresponding
transfer functions coincide. Thus we have that $\left[
                                                                                      \begin{array}{cc}
                                                                                        \vpj & 0 \\
                                                                                      \end{array}
                                                                                    \right]
$ coincides with $\sjp$. As a consequence, there exist unitary
operators $\be : \clY \rightarrow \clW$ and $\hat{\al} : \clN \op
\clH \rightarrow \mathcal{N}^\prime$ such that $\hat{\al}
(~\clU^{0 \perp}_\vp) = \clV^{0 \perp}_\psi$, satisfying the
coincidence identity. Since $\mathrm{dim}~\clU^0_\vp = \mathrm{dim}~\clV^0_\psi$,
choose any unitary $\tau : \clU^0_\vp \rightarrow \clV^0_\psi$ and
consider the unitary operator
$$\al = \left(
          \begin{array}{cc}
            \hat{\al}|_{\clU^{0 \perp}_\vp} & 0 \\
            0 & \tau \\
          \end{array}
        \right) : \clU^{0 \perp}_\vp \op \clU^0_\vp \rightarrow \clV~.
$$

Finally, the coincidence of $\vp$ and $\psi$ follow from the easy
computation below.

Let $h\in \clU$ and $z\in\Bd$, then

\begin{eqnarray*}
\psi(z) \al~ h &=& \psi(z) \al (P_{\clU^{0 \perp}_\vp} h \op P_{\clU^{0}_\vp}h)\\
&=& \psi(z) (\hat{\al}|_{\clU^{0 \perp}_\vp} P_{\clU^{0 \perp}_\vp} h \op \tau P_{\clU^{0}_\vp}h )\\
&=& \psi(z) (\hat{\al}|_{\clU^{0 \perp}_\vp} P_{\clU^{0 \perp}_\vp} h) \\
&=& \psi(z) J^\prime \hat{\al} (P_{\clU^{0 \perp}_\vp} h \op 0)\\
&=& \sjp(z) \hat{\al} (P_{\clU^{0 \perp}_\vp} h \op 0) \\
&=& \be [~~\vpj~~~0~~] (P_{\clU^{0 \perp}_\vp} h \op 0)\\
&=& \be \vpj(z)  (P_{\clU^{0 \perp}_\vp} h)\\
&=& \be \vp(z) (P_{\clU^{0 \perp}_\vp} h) = \be \vp(z)~ h~.
\end{eqnarray*}

That proves $(2) \Rightarrow (1)$. The converse proof will follow
from the following theorem about coincidence of Schur multipliers
admitting coisometric functional model realizations. When the
multipliers admit coisometric functional model realizations, the
partial isometry $J$ does not figure in the characterization.

\bt \label{coiso}

Let $\vp\in\sduy$ and $\psi\in\sdvw$ be such that both admit coisometric functional model realizations. Then the following are equivalent :
\begin{enumerate}
  \item $\vp$ coincides with $\psi$.
  \item $\mathrm{dim}~\clU^0_\vp = \mathrm{dim}~\clV^0_\psi$ and for any two arbitrary choices of coisometric functional model realizations say $U^\vp$ and $U^\psi$ of $\vp$ and $\psi$ respectively, there exists an auxiliary Hilbert space $\clM$ such that either $U^\vp_\clM$ is unitarily coincident to $U^\psi$ or $U^\psi_\clM$ is unitarily coincident to $U^\vp$.
\end{enumerate}

\et

\vspace{3mm}

Proof : $(1) \Rightarrow (2)$. Let $\vp$ and $\psi$ be coincident. Then there exist unitaries say, $\al : \clU \rightarrow \clV$ and $\be : \clY \rightarrow \clW$, such that $$\be \vp(z)=\psi(z)\al$$holds for all $z\in\Bd$.

~~~~~

Let $U^\vp$ and $U^\psi$ be any two coisometric functional model realizations of $\vp$ and $\psi$ respectively. By Theorem 2.7 of \cite{BBF1}, the colligation $U^\vp$ can be explicitly described as :
$$U^\vp = \left(
            \begin{array}{cc}
              A^\vp & B^\vp \\
              C^\vp & \vp(0) \\
            \end{array}
          \right) : \left(
                      \begin{array}{c}
                        \hkp \\
                        \clU \\
                      \end{array}
                    \right) \rightarrow \left(
                                          \begin{array}{c}
                                            \hkp^d \\
                                            \clY \\
                                          \end{array}
                                        \right)
$$where the non-unique operator $B^{\vp *} : \hkp^d \rightarrow \clU$ is given by
$$\left(
    \begin{array}{cc}
      -G_2^\vp T_{12}^{\vp *} G_1^{\vp *} + \zeta^\vp (I_{\Dp^\perp} - G_1^{\vp} G_1^{\vp *} )^{\frac{1}{2}} & R^\vp \\
    \end{array}
  \right) : \Dp^\perp \op \Dp \rightarrow \clU,
$$for some isometry $\zeta^\vp : \overline{\mathrm{ran}}~(I_{\Dp^\perp} - G_1^{\vp} G_1^{\vp *} )^{\frac{1}{2}} \rightarrow \clU^0_\vp$. The non-uniqueness of $B^\vp$ comes from dependence of the parameter $\zeta^\vp$.

~~~~~

Similarly, one has a description of the chosen colligation $U^\psi$ in terms of the operators $G_2^\psi, T_{12}^{\psi *}, G_1^{\psi}$ and some isometric operator $\zeta^\psi : \overline{\mathrm{ran}}~(I_{\Ds^\perp} - G_1^{\psi} G_1^{\psi *} )^{\frac{1}{2}} \rightarrow \clV^0_\vp$.

~~~~~

Due to coincidence of the Schur multipliers one has
\begin{eqnarray*}
 \Pd|_{\Dp^\perp}~(I_{\Dp^\perp} - G_1^{\vp} G_1^{\vp *} )^{\frac{1}{2}} &=& (I_{\Ds^\perp} - G_1^{\psi} G_1^{\psi *} )^{\frac{1}{2}}~\Pd|_{\Dp^\perp}~.
\end{eqnarray*}

Thus,
$$\Pd|_{\Dp^\perp} (~\overline{\mathrm{ran}}~(I_{\Dp^\perp} - G_1^{\vp} G_1^{\vp *}~)^{\frac{1}{2}}) = \overline{\mathrm{ran}}~(I_{\Ds^\perp} - G_1^{\psi} G_1^{\psi *} )^{\frac{1}{2}}~.$$

~~~~~

If
\begin{equation} \label{dim}
\mathrm{dim}~(\clU^0_\vp \ominus \mathrm{ran}~\zeta^\vp) \leq \mathrm{dim}~(\clV^0_\psi \ominus \mathrm{ran}~\zeta^\psi)~,
\end{equation}
then choose a Hilbert space $\clM$ such that
$$\mathrm{dim}~[\clM \op (\clU^0_\vp \ominus \mathrm{ran}~\zeta^\vp)] = \mathrm{dim}~(\clV^0_\psi \ominus \mathrm{ran}~\zeta^\psi)~.$$

~~~~~

By virtue of the fact that $\zeta^{(\cdot)}$ is an isometry (and hence a unitary onto its range) choose a unitary $\delta : \mathrm{ran}~\zeta^\vp \rightarrow \mathrm{ran}~\zeta^\psi$ and $\delta^\prime : \clM \op (\clU^0_\vp \ominus \mathrm{ran}~\zeta^\vp) \rightarrow \clV^0_\psi \ominus \mathrm{ran}~\zeta^\psi$ such that $\delta \op \delta^\prime : \clU^0_{\vp} \op \clM \rightarrow \clV^0_\psi$ is unitary and the following diagram commutes.

\setlength{\unitlength}{3mm}
 \begin{center}
 \begin{picture}(20,14)(0,0)
 \put(-6.5,2.5){$\overline{\mathrm{ran}}~(I_{\mathcal{D}_{\psi}^\perp} - G_1^\psi G_1^{\psi *})^{\frac{1}{2}}$} \put(13.5,2.5){$\mathcal{V}^0_\psi$}
 \put(8.6,1.6){$\zeta^\psi$}
 \put(-4.0,6.5){$ \Gamma^d|_{\mathcal{D}^\perp_\varphi}$} \put(14.5,6.5){$\delta \oplus \delta^\prime$}
 \put(-6.5,10){$\overline{\mathrm{ran}}~(I_{\mathcal{D}_{\varphi}^\perp} - G_1^\varphi G_1^{\varphi *})^{\frac{1}{2}}$} \put(13.5,10){$\mathcal{U}^0_\varphi \oplus \mathcal{M}$}
 \put(8.6,11){$\zeta^\varphi$}
 \put(5.5,3.0){ \vector(1,0){7}} \put(5.5,10.5){ \vector(1,0){7}}
 \put(-0.8,9.2){ \vector(0,-1){5}} \put(13.5,9.2){ \vector(0,-1){5}}

 \end{picture}
 \end{center}

As $\al (~\clU^{0 \perp}_\vp) = \clV^{0 \perp}_\psi$, we consider the unitary operator $\tilde{\al} : \clU \op \clM \rightarrow \clV$ given by
$$\tilde{\al} = \left(
                  \begin{array}{cc}
                    \al|_{\clU^{0 \perp}_\vp} & 0 \\
                    0 & \delta \op \delta^\prime \\
                  \end{array}
                \right) : \clU \op \clM \rightarrow \clV~.
$$

~~~~~

By definition of $\tilde{\al}$ we have
$$\be \vp_\clM (z) = \psi (z) \tilde{\al}$$for all $z\in\Bd$.

~~~~~

Consider the operator $\left[
                         \begin{array}{cc}
                           B^\vp & 0 \\
                         \end{array}
                       \right]^* : \hkp^d \rightarrow \clU \op \clM
$ acting as :
$$h \longmapsto B^{\vp *} h \op 0~.$$

~~~~~

For $h=h_1 \op h_2$ where $h_1 \in \Dp$ and $h_2 \in \Dp^\perp$ one has
$$\left[
     \begin{array}{cc}
       B^\vp & 0 \\
     \end{array}
   \right]^* h = (-G_2^\vp T_{12}^{\vp *} G_1^{\vp *} h_2 + R^\vp h_1) \op (\zeta^\vp (I_{\Dp^\perp} - G_1^{\vp} G_1^{\vp *} )^{\frac{1}{2}} h_2 + 0)~.$$
The first term in brackets belongs to $\clU^{0 \perp}_\vp$ and the second term belongs to $\clU^0_\vp \op \clM$.

~~~~~

Taking into account the fact (due to coincidence of $\vp$ and $\psi$) that
$$- \al~ G_2^\vp T_{12}^{\vp *} G_1^{\vp *} = -G_2^\psi T_{12}^{\psi *} G_1^{\psi *}~ \Pd|_{\Dp^\perp}$$
and from the construction of the unitary $\tilde{\al}$ it is easy to see that
$$\tilde{\al} \left[
                 \begin{array}{cc}
                   B^\vp & 0 \\
                 \end{array}
               \right]^* = B^{\psi *}~ \Pd~.$$

~~~~~

These relations lead to the fact that
$$\left(
    \begin{array}{cc}
      \Pd & 0 \\
      0 & \be \\
    \end{array}
  \right) U^\vp_\clM = U^\psi \left(
                              \begin{array}{cc}
                                \Gamma & 0 \\
                                0 & \tilde{\al} \\
                              \end{array}
                            \right)
$$and our claim is proved. If, in (\ref{dim}), the inequality holds the other way, then the same argument with $\vp$ and $\psi$ interchanged, would give us the unitary coincidence of $U^\psi_{\clM}$ and $U^\vp$.

~~~~~

$(2) \Rightarrow (1)$. Without loss of generality, let us assume that there exists a Hilbert space $\clM$ such that $U^\vp_\clM$ is unitarily coincident to $U^\psi$. Then we know that the corresponding transfer functions coincide. But the transfer function of $U^\vp_\clM$ is $\vp_\clM$.

~~~~~

Let $\hat{\al} : \clU \op \clM \rightarrow \clV$ and $\be : \clY \rightarrow \clW$ be unitaries such that
$$\be \vp_\clM (z)=\psi (z) \hat{\al}$$holds for all $z\in\Bd$.

~~~~~

From above we conclude that $\hat{\al} (~\clU^{0 \perp }_\vp) = \clV^{0 \perp}_\psi$, since $(\clU \op \clM)^0_{\vp_{\clM}} = \clU^0_\vp \op \clM$. Also, by virtue of the assumption $\mathrm{dim}~\clU^0_\vp = \mathrm{dim}~\clV^0_\psi$ choose a unitary $\tau : \clU^0_\vp \rightarrow \clV^0_\psi$ and consider the unitary operator
$$\al = \left(
          \begin{array}{cc}
            \hat{\al}|_{\clU^{0 \perp }_\vp} & 0 \\
            0 & \tau \\
          \end{array}
        \right) : \clU \rightarrow \clV~.
$$

Now it is easy to see that $\be \vp (z)=\psi (z)\al$ holds for all
$z\in\Bd$ and our claim is proved. \qed

Now we complete the proof of the main theorem.

Let $\vp \in \sduy$ and $\psi \in \sdvw$ and let $\vp$ coincide
with $\psi$. Let $\al$ and $\be$ the unitary operators as before.
Also let $\clN^\prime$ and $J^\prime$ be the Hilbert space and the
partial isometry counterparts of $\clN$ and $J$ respectively for
the Schur multiplier $\psi$. Then it is easy to see that the
operator
$$\al \op \Pd|_{\Dp^\perp}$$is a unitary operator from $\clN$ onto $\clN^\prime$.

Considering the above unitary operator it is easy to see that
$\vpj$ and $\sjp$ coincide. Indeed, for $m\in\clN$ and $z\in\Bd$,

\begin{eqnarray*}
\sjp (z) (\al \op \Pd|_{\Dp^\perp}) m &=& \psi (z) J^\prime (\al \op \Pd|_{\Dp^\perp}) m \\
&=& \psi (z) \al (P_{\clU^{0 \perp}_\vp}~ m) \\
&=& \be \vp (z) (P_{\clU^{0 \perp}_\vp}~ m) \\
&=& \be \vp (z) J m \\
&=& \be \vpj (z) m~.
\end{eqnarray*}

Now to complete the proof  $(1) \Rightarrow (2)$ of the main
theorem, if $(1)$ holds then $\vpj$ coincides with $\sjp$.
Considering their coisometric functional model realizations
$U^{\vpj}$ and $U^{\sjp}$ obtained canonically from an arbitary
but fixed choice $U^\vp$ and $U^\psi$ of functional model
realizations of $\vp$ and $\psi$, an application of Theorem
\ref{coiso} proves our claim. \qed

\vspace{1cm}

\textbf{Acknowledgement.} The second named author acknowledges Ramanna
Fellowship from DST and UGC SAP Phase IV.

\end{document}